\ifx\shlhetal\undefinedcontrolsequence\let\shlhetal\relax\fi
\documentclass{amsart}
\usepackage{amsmath} 
\usepackage{amssymb}

\newtheorem{theorem}{Theorem}
\newtheorem{claim}{Claim}[theorem]
\newtheorem{lemma}[theorem]{Lemma} 
\newtheorem{proposition}[theorem]{Proposition}

\theoremstyle{definition}
\newtheorem{definition}[theorem]{Definition}

\newtheorem{problem}[theorem]{Problem}

\theoremstyle{remark}

\newtheorem{conclusion}[theorem]{Conclusion}

\newcommand{\bB}{{\mathbb B}}
\newcommand{\bP}{{\mathbb P}}

\newcommand{\bQ}{{\mathbb Q}}

\newcommand{\comp}{\circ}
\newcommand{\forces}{\Vdash}
\newcommand{\bV}{{\bf V}}
\newcommand{\rest}{\restriction}
\newcommand{\dom}{{\rm dom}}

\newcommand{\otp}{{\rm otp}}
\newcommand{\cl}{{\rm cl}}
\newcommand{\vare}{\varepsilon}
\newcommand{\dbtl}{\dot{\bB}^\theta_\lambda}
\newcommand{\Dep}{{\rm Depth}}
\newcommand{\rht}{{\rm ht}}
\newcommand{\ptl}{{{\mathbb P}^\theta_\lambda}}

\title{Historic forcing for $\Dep$}
\author{Andrzej Ros{\l}anowski}
\address{Department of Mathematics\\
 University of Nebraska at Omaha\\
 Omaha, NE 68182-0243, USA\\
 and Mathematical Institute of Wroclaw University\\
 50384 Wroclaw, Poland} 
\email{roslanowski@unomaha.edu}
\urladdr{http://www.unomaha.edu/$\sim$aroslano}
\thanks{The first author thanks the KBN (Polish Committee of  Scientific
 Research) for partial support through grant 2 P03 A 01109.} 

\author{Saharon Shelah}
\address{Institute of Mathematics\\
 The Hebrew University of Jerusalem\\
 91904 Jerusalem, Israel\\
 and  Department of Mathematics\\
 Rutgers University\\
 New Brunswick, NJ 08854, USA}
\email{shelah@math.huji.ac.il}
\urladdr{http://www.math.rutgers.edu/$\sim$shelah}
\thanks{The research of the second author was partially supported by the
 Israel Science Foundation. Publication 733} 

\subjclass{Primary 03E35, 03G05; Secondary 03E05, 06Exx}
\keywords{Boolean algebras, depth, historic forcing}

\begin{document}

\begin{abstract}
We show that, consistently, for some regular cardinals $\theta<\lambda$,
there exist a Boolean algebra $\bB$ such that $|\bB|=\lambda^+$ and for
every subalgebra $\bB'\subseteq\bB$ of size $\lambda^+$ we have $\Dep(\bB')=
\theta$. 
\end{abstract}

\maketitle

\section{Introduction}
The present paper is concerned with forcing a Boolean algebra which has some
prescribed properties of $\Dep$. Let us recall that, for a Boolean algebra
$\bB$, its depth is defined as follows:
\[\begin{array}{lcl}
\Dep(\bB)&=&\sup\{|X|: X\subseteq\bB\mbox{ is well-ordered by the Boolean
ordering}\;\},\\ 
\Dep^+(\bB)&=&\sup\{|X|^+: X\subseteq\bB\mbox{ is well-ordered by the
Boolean ordering}\;\}. 
  \end{array}\]
($\Dep^+(\bB)$ is used to deal with attainment properties in the definition
of $\Dep(\bB)$, see e.g.~\cite[\S 1]{RoSh:534}.) The depth (of Boolean
algebras) is among cardinal functions that have more algebraic origins, and
their relations to ``topological fellows'' is often indirect, though
sometimes very surprising. For example, if we define
\[\Dep_{{\rm H}+}(\bB)=\sup\{\Dep(\bB/I): I\mbox{ is an ideal in }\bB\;\},\]
then for any (infinite) Boolean algebra $\bB$ we will have that $\Dep_{{\rm
H}+}(\bB)$ is the tightness $t(\bB)$ of the algebra $\bB$ (or the tightness
of the topological space ${\rm Ult}(\bB)$ of ultrafilters on $\bB$), see
\cite[Theorem 4.21]{M2}. A somewhat similar function to $\Dep_{{\rm H}+}$
is obtained by taking $\sup\{\Dep(\bB'): \bB'$ is a subalgebra of $\bB\;\}$,
but clearly this brings nothing new: it is the old Depth. But if one wants
to understand the behaviour of the depth for subalgebras of the considered
Boolean algebra, then looking at the following {\em subalgebra $\Dep$
relation\/} may be very appropriate:
\[\begin{array}{lr}
\Dep_{\rm Sr}(\bB)=\{(\kappa,\mu):&\mbox{there is an infinite subalgebra
$\bB'$ of $\bB$ such that }\ \\
&|\bB'|=\mu\mbox{ and }\Dep(\bB')=\kappa\;\}.
  \end{array}\]
A number of results related to this relation is presented by Monk in 
\cite[Chapter 4]{M2}. There he asks if there are a Boolean algebra $\bB$ and
an infinite cardinal $\theta$ such that $(\theta,(2^\theta)^+)\in \Dep_{\rm
Sr}(\bB)$, while $(\omega, (2^\theta)^+)\notin\Dep_{\rm Sr}(\bB)$ (see Monk
\cite[Problem 14]{M2}; we refer the reader to Chapter 4 of Monk's book
\cite{M2} for the motivation and background of this problem). Here we will
partially answer this question, showing that it is consistent that there is
such $\bB$ and $\theta$. The question if that can be done in ZFC remains
open. 

Our consistency result is obtained by forcing, and the construction of the
required forcing notion is interesting {\em per se}. We use the method of
{\em historic forcing\/} which was first applied  in Shelah and Stanley
\cite{ShSt:258}. The reader familiar with \cite{ShSt:258} will notice
several correspondences between the construction here and the method
used there. However, we do not relay on that paper and our presentation
here is self-contained. 

Let us describe how our historic forcing notion is built. So, we fix two
(regular) cardinals $\theta,\lambda$ and our aim is to force a Boolean
algebra $\dbtl$ such that $|\dbtl|=\lambda^+$ and for every subalgebra
$\bB\subseteq\dbtl$ of size $\lambda^+$ we have $\Dep(\bB)=\theta$. The
algebra $\dbtl$ will be generated by $\langle x_i:i\in\dot{U}\rangle$ for
some set $\dot{U}\subseteq\lambda^+$. A condition $p$ will be an
approximation to the algebra $\dbtl$, it will carry the information on what
is the subalgebra $\bB_p=\langle x_i: i\in u^p\rangle_{\dbtl}$ for some
$u^p\subseteq\lambda^+$. A natural way to describe algebras in this context
is by listing ultrafilters (or: homomorphisms into $\{0,1\}$):
\begin{definition}
\label{0.C}
For a set $w$ and a family $F\subseteq 2^{\textstyle w}$ we define

\noindent $\cl(F)=\{g\in 2^{\textstyle w}: (\forall u\in [w]^{\textstyle
<\omega})(\exists f\in F)(f\rest u=g\rest u)\}$,

\noindent $\bB_{(w,F)}$ is the Boolean algebra generated freely by
$\{x_\alpha:\alpha\in w\}$ except that

if $u_0,u_1\in [w]^{\textstyle <\omega}$ and there is no $f\in F$ such
that $f\rest u_0\equiv 0$, $f\rest u_1\equiv 1$

then $\bigwedge\limits_{\alpha\in u_1} x_\alpha\wedge 
\bigwedge\limits_{\alpha\in u_0} (-x_\alpha)=0$.
\end{definition}
This description of algebras is easy to handle, for example: 
\begin{proposition}
[see {\cite[2.6]{Sh:479}}]
\label{0.D}
Let $F\subseteq 2^{\textstyle w}$. Then:
\begin{enumerate}
\item Each $f\in F$ extends (uniquely) to a homomorphism from $\bB_{(w,F)}$
to $\{0,1\}$ (i.e.~it preserves the equalities from the definition of
$\bB_{(w,F)}$). If $F$ is closed, then every homomorphism from $\bB_{(w,F)}$
to $\{0,1\}$ extends exactly one element of $F$.  
\item If $\tau(y_0,\ldots,y_\ell)$ is a Boolean term and $\alpha_0,\ldots,
\alpha_\ell\in w$ are distinct then
\[\begin{array}{l}
\bB_{(w,F)}\models\tau(x_{\alpha_0},\ldots,x_{\alpha_\ell})\neq 0\qquad
\qquad\mbox{ if and only if}\\
(\exists f\in F)(\{0,1\}\models\tau(f(\alpha_0),\ldots,f(\alpha_k))=1).
  \end{array}\]
\item If $w\subseteq w^*$, $F^*\subseteq 2^{\textstyle w^*}$ and
\[(\forall f\in F)(\exists g\in F^*)(f\subseteq g)\quad\mbox{ and }\quad
(\forall g\in F^*)(g\rest w\in\cl(F))\]
then $\bB_{(w,F)}$ is a subalgebra of $\bB_{(w^*,F^*)}$.
\end{enumerate}
\end{proposition}
So each condition $p$ in our forcing notion $\ptl$ will have a set $u^p\in
[\lambda^+]^{<\lambda}$ and a closed set $F^p\subseteq 2^{\textstyle u^p}$
(and the respective algebra will be $\bB_p=\bB_{(u^p,F^p)}$). But to make
the forcing notion work, we will have to put more restrictions on our
conditions, and we will be taking only those conditions that have to be
taken to make the arguments work. For example, we want that cardinals are
not collapsed by our forcing, and demanding that $\ptl$ is $\lambda^+$-cc
(and somewhat $({<}\lambda)$--closed) is natural in this context. How do we
argue that a forcing notion is $\lambda^+$--cc? Typically we start with a
sequence of $\lambda^+$ distinct conditions, we carry out some ``cleaning
procedure'' (usually involving the $\Delta$--lemma etc), and we end up with
(at least two) conditions that ``can be put together''. Putting together two
(or more) conditions that are approximations to a Boolean algebra means
amalgamating them. There are various ways to amalgamate conditions - we will
pick one that will work for several purposes. Then, once we declare that
some conditions forming a ``clean'' $\Delta$--sequence of length $\theta$
are in $\ptl$, we will be bound to declare that the amalgamation is in our
forcing notion. The amalgamation (and natural limits) will be the only way
to build new conditions from the old ones, but the description above still
misses an important factor. So far, a condition does not have to know what
are the reasons for it to be called to $\ptl$. This information is {\em the
history of the condition\/} and it will be encoded by two functions
$h^p,g^p$. (Actually, these functions will give histories of all elements of
$u^p$ describing why and how those points were incorporated to $u^p$. Thus
both functions will be defined on $u^p\times\rht(p)$, were $\rht(p)$ is the
height of the condition $p$, that is the step in our construction at which
the condition $p$ is created.) We will also want that our forcing is
suitably closed, and getting ``$({<}\lambda)$--strategically closed'' would
be fine. To make that happen we will have to deal with two relations on on
$\ptl$: $\leq_{\rm pr}$ and $\leq$. The first (``pure'') is
$({<}\lambda)$--closed and it will help in getting the strategic closure of
the second (main) one. In some sense, the relation $\leq_{\rm pr}$
represents ``the official line in history'', and sometimes we will have to
rewrite that official history, see Definition \ref{defptran} and Lemma
\ref{3y1} (on changing history see also Orwell \cite{Or49}). 

The forcing notion $\ptl$ has some other interesting features. (For example, 
conditions are very much like fractals, they contain many self-similar pieces
(see Definition \ref{defcompo} and Lemma \ref{3.4x}).) The method of
historic forcing notions could be applicable to more problems, and this is
why in our presentation we separated several observations of general
character (presented in the first section) from the problem specific
arguments (section 2)
\medskip

\noindent{\bf Notation:}\qquad Our notation is standard and compatible with
that of classical textbooks on set theory (like Jech \cite{J}) and Boolean
algebras (like Monk \cite{M1}, \cite{M2}). However in forcing considerations
we keep the older tradition that
\begin{center}
{\em
the stronger condition is the greater one.
}
\end{center}
Let us list some of our notation and conventions.

\begin{enumerate}
\item Throughout the paper, $\theta,\lambda$ are fixed regular infinite
cardinals, $\theta<\lambda$.  
\item A name for an object in a forcing extension is denoted with a dot
above (like $\dot{X}$) with one exception: the canonical name for a generic
filter in a forcing notion $\bP$ will be called $\Gamma_\bP$. For a
$\bP$--name $\dot{X}$ and a $\bP$--generic filter $G$ over $\bV$, the
interpretation of the name $\dot{X}$ by $G$ is denoted by $\dot{X}^G$. 
\item $i,j,\alpha,\beta,\gamma,\delta,\ldots$ will denote ordinals.
\item For a set $X$ and a cardinal $\lambda$, $[X]^{\textstyle<\lambda}$
stands for the family of all subsets of $X$ of size less than $\lambda$. The
family of all functions from $Y$ to $X$ is called $X^{\textstyle Y}$. If 
$X$ is a set of ordinals then its order type is denoted by $\otp(X)$.
\item In Boolean algebras we use $\vee$ (and $\bigvee$), $\wedge$ (and
$\bigwedge$) and $-$ for the Boolean operations. If $\bB$ is a Boolean
algebra, $x\in\bB$ then $x^0=x$, $x^1=-x$. 
\item For a subset $Y$ of an algebra $\bB$, the subalgebra of $\bB$
generated by $Y$ is denoted by $\langle Y\rangle_{\bB}$.
\end{enumerate}

\noindent{\bf Acknowledgements:}\qquad We would like to thank the referee for
valuable comments and suggestions. 

\section{The forcing and its basic properties}
Let us start with the definition of the forcing notion $\ptl$. By induction
on $\alpha<\lambda$ we will define sets of conditions $P^{\theta,
\lambda}_\alpha$, and for each $p\in P^{\theta,\lambda}_\alpha$ we will
define $u^p,F^p,\rht(p),h^p$ and $g^p$. Also we will define relations
$\leq^\alpha$ and $\leq^\alpha_{\rm pr}$ on $P^{\theta,\lambda}_\alpha$. Our
inductive requirements are: 
\begin{enumerate}
\item[(i)$_\alpha$] for each $p\in P^{\theta,\lambda}_\alpha$:\\
$u^p\in [\lambda^+]^{\textstyle<\lambda}$, $\rht(p)\leq\alpha$, $F^p
\subseteq 2^{\textstyle u^p}$ is a non-empty closed set, $g^p$ is a function
with domain $\dom(g^p)=u^p\times\rht(p)$ and values of the form $(\ell,
\tau)$, where $\ell<2$ and $\tau$ is a Boolean term, and $h^p:u^p\times
\rht(p)\longrightarrow\theta+2$ is a function,
\item[(ii)$_\alpha$] $\leq^\alpha,\leq^\alpha_{\rm pr}$ are transitive and
reflexive relations on $P^{\theta,\lambda}_\alpha$, and $\leq^\alpha$
extends $\leq^\alpha_{\rm pr}$,
\item[(iii)$_\alpha$] if $p,q\in P^{\theta,\lambda}_\alpha$, $p\leq^\alpha
q$, then $u^p\subseteq u^q$, $\rht(p)\leq\rht(q)$, and $F^p=\{f\rest u^p:
f\in F^q\}$, and if $p\leq^\alpha_{\rm pr} q$, then for every $i\in u^p$ and
$\xi<\rht(p)$ we have $h^p(i,\xi)=h^q(i,\xi)$ and $g^p(i,\xi)= g^q(i,\xi)$, 
\item[(iv)$_\alpha$] if $\beta<\alpha$ then $P^{\theta,\lambda}_\beta
\subseteq P^{\theta,\lambda}_\alpha$, and $\leq^\alpha_{\rm pr}$ extends
$\leq^\beta_{\rm pr}$, and $\leq^\alpha$ extends $\leq^\beta$.  
\end{enumerate}
For a condition $p\in P^{\theta,\lambda}_\alpha$, we will also declare
that $\bB^p=\bB_{(u^p,F^p)}$ (the Boolean algebra defined in Definition
\ref{0.C}).    
\medskip

We define $P^{\theta,\lambda}_0=\{\langle\xi\rangle:\xi<\lambda^+\}$ and
for $p=\langle\xi\rangle$ we let $F^p=2^{\textstyle\{\xi\}}$, $\rht(p)=0$
and $h^p=\emptyset=g^p$. The relations $\leq^0_{\rm pr}$ and $\leq^0$ both
are the equality. [Clearly these objects are as declared, i.e, clauses
(i)$_0$--(iv)$_0$ hold true.]    
\medskip

If $\gamma<\lambda$ is a limit ordinal, then we put 
\[\begin{array}{l}
P^*_\gamma=\big\{\langle p_\xi:\xi<\gamma\rangle: (\forall\xi<\zeta<
\gamma)(p_\xi\in P^{\theta,\lambda}_\xi\ \&\ \rht(p_\xi)=\xi \ \&\ p_\xi
\leq^\zeta_{\rm pr} p_\zeta)\big\},\\
P^{\theta,\lambda}_\gamma=\bigcup\limits_{\alpha<\gamma}P^{\theta,
\lambda}_\alpha\cup P^*_\gamma,
  \end{array}\]
and for $p=\langle p_\xi:\xi<\gamma\rangle\in P^*_\gamma$ we let
\[u^p=\bigcup\limits_{\xi<\gamma} u^{p_\xi},\quad F^p=\{f\in 2^{\textstyle
u^p}:(\forall\xi<\gamma)(f\rest u^{p_\xi}\in F^{p_\xi})\},\quad \rht(p)=
\gamma\]
and $h^p=\bigcup\limits_{\xi<\gamma} h^{p_\xi}$ and $g^p=\bigcup\limits_{\xi<
\gamma} g^{p_\xi}$. We define $\leq^\gamma$ and $\leq^\gamma_{\rm pr}$ by: 
\medskip

\noindent $p\leq^\gamma_{\rm pr}q$\quad if and only if

{\em either\/} $p,q\in P^{\theta,\lambda}_\alpha$, $\alpha<\gamma$ and
$p\leq^\alpha_{\rm pr}q$, 

{\em or\/} $q=\langle q_\xi:\xi<\gamma\rangle\in P^*_\gamma$, $p\in
P^{\theta,\lambda}_\alpha$ and $p\leq^\alpha_{\rm pr}q_\alpha$ for some
$\alpha<\gamma$, 

{\em or\/} $p=q$;
\smallskip

\noindent $p\leq^\gamma q$\quad if and only if

{\em either\/} $p,q\in P^{\theta,\lambda}_\alpha$, $\alpha<\gamma$ and
$p\leq^\alpha q$, 

{\em or\/} $q=\langle q_\xi:\xi<\gamma\rangle\in P^*_\gamma$, $p\in
P^{\theta,\lambda}_\alpha$ and $p\leq^\alpha q_\alpha$ for some $\alpha<
\gamma$,  

{\em or\/} $p=\langle p_\xi:\xi<\gamma\rangle\in P^*_\gamma$, $q=\langle
q_\xi:\xi<\gamma\rangle\in P^*_\gamma$ and 
\[(\exists\delta<\gamma)(\forall\xi<\gamma)(\delta\leq\xi\ \Rightarrow\
p_\xi\leq^\xi q_\xi).\]

\noindent [It is straightforward to show that clauses
(i)$_\gamma$--(iv)$_\gamma$ hold true.]
\medskip

Suppose now that $\alpha<\lambda$. Let $P^*_{\alpha+1}$ consist of all
tuples  
\[\langle\zeta^*,\tau^*,n^*,u^*,\langle p_\xi,v_\xi:\xi<\theta\rangle\rangle\]
such that for each $\xi_0<\xi_1<\theta$:
\begin{enumerate}
\item[$(\alpha)$] $\zeta^*<\theta$, $n^*<\omega$, $\tau^*=\tau^*(y_1,\ldots,y_{
n^*})$ is a Boolean term, $u^*\in [\lambda^+]^{\textstyle<\lambda}$,
\item[$(\beta)$]  $p_{\xi_0}\in P^{\theta,\lambda}_\alpha$, $\rht(p)=\alpha$,
$v_{\xi_0}\in [u^{p_{\xi_0}}]^{\textstyle n^*}$,
\item[$(\gamma)$] the family $\{u^{p_\xi}:\xi<\theta\}$ forms a
$\Delta$--system with heart $u^*$ and $u^{p_{\xi_0}}\setminus u^*\neq
\emptyset$ and
\[\sup(u^*)<\min(u^{p_{\xi_0}}\setminus u^*)\leq\sup (u^{p_{\xi_0}}\setminus
u^*)<\min(u^{p_{\xi_1}}\setminus u^*),\]
\item[$(\delta)$] $\otp(u^{p_{\xi_0}})=\otp(u^{p_{\xi_1}})$ and if $H:u^{p_{
\xi_0}}\longrightarrow u^{p_{\xi_1}}$ is the order isomorphism then $H\rest
u^*$ is the identity on $u^*$, $F^{p_{\xi_0}}=\{f\comp H:f\in F^{p_{\xi_1}}
\}$, $H[v_{\xi_0}]=v_{\xi_1}$ and  
\[(\forall j\in u^{p_{\xi_0}})(\forall\beta<\alpha)(h^{p_{\xi_0}}(j,\beta)=
h^{p_{\xi_1}}(H(j),\beta)\ \&\ g^{p_{\xi_0}}(j,\beta)=g^{p_{\xi_1}}(H(j),
\beta)).\]
\end{enumerate}
We put $P^{\theta,\lambda}_{\alpha+1}=P^{\theta,\lambda}_\alpha\cup P^*_{
\alpha+1}$ and for $p=\langle\zeta^*,\tau^*,n^*,u^*,\langle p_\xi,v_\xi:\xi<
\theta\rangle\rangle\in P^*_{\alpha+1}$ we let $u^p=\bigcup\limits_{\xi<
\theta} u^{p_\xi}$ and
\[\begin{array}{ll}
F^p=\{f\in 2^{\textstyle u^p}:& (\forall\xi<\theta)(f\rest u^{p_\xi}\in
F^{p_{\xi}})\mbox{ and for all }\xi<\zeta<\theta\\
\ &f(\sigma_{\rm maj}(\tau_{3\cdot\xi},\tau_{3\cdot\xi+1},\tau_{3\cdot\xi+
2}))\leq f(\sigma_{\rm maj}(\tau_{3\cdot\zeta},\tau_{3\cdot\zeta+1},\tau_{3
\cdot\zeta+2}))\},
  \end{array}\]
where $\tau_\xi=\tau^*(x_i:i\in v_\xi)$ for $\xi<\theta$ (so $\tau_\xi$ is
an element of the algebra $\bB^{p_\xi}=\bB_{(u^{p_\xi},F^{p_\xi})}$), and
$\sigma_{\rm maj}(y_0,y_1,y_2)=(y_0\wedge y_1)\vee (y_0\wedge y_2)\vee
(y_1\wedge y_2)$. Next we let $\rht(p)=\alpha+1$ and we define functions
$h^p,g^p$ on $u^p\times (\alpha+1)$ by
\[h^p(j,\beta)=\left\{\begin{array}{lll}
h^{p_{\xi}}(j,\beta)&\mbox{if}&j\in u^{p_\xi},\ \xi<\theta,\ \beta<\alpha,\\
\theta                 &\mbox{if}&j\in u^*,\ \beta=\alpha,\\
\theta+1               &\mbox{if}&j\in u^{p_{\zeta^*}}\setminus u^*,\
\beta=\alpha,\\ 
\xi                    &\mbox{if}&j\in u^{p_\xi}\setminus u^*,\ \xi<\theta,\
\xi\neq\zeta^*,\ \beta=\alpha, 
			 \end{array}\right.\]
\[g^p(j,\beta)=\left\{\begin{array}{lll}
g^{p_{\xi}}(j,\beta)&\mbox{if}&j\in u^{p_\xi},\ \xi<\theta,\ \beta<\alpha,\\
(1,\tau^*)          &\mbox{if}&j\in v_\xi,\ \xi<\theta,\ \beta=\alpha,\\
(0,\tau^*)          &\mbox{if}&j\in u^{p_\xi}\setminus v_\xi,\ \xi<\theta,\
\beta=\alpha.
			 \end{array}\right.\]
Next we define the relations $\leq^{\alpha+1}_{\rm pr}$ and $\leq^{\alpha+1}$
by: 
\medskip

\noindent $p\leq^{\alpha+1}_{\rm pr}q$\quad if and only if

{\em either\/} $p,q\in P^{\theta,\lambda}_\alpha$ and $p\leq^\alpha_{\rm
pr}q$, 

{\em or\/} $q=\langle \zeta^*,\tau^*,n^*,u^*,\langle q_\xi,v_\xi:\xi<\theta
\rangle\rangle\in P^*_{\alpha+1}$, $p\in P^{\theta,\lambda}_\alpha$, and
$p\leq^\alpha_{\rm pr} q_{\zeta^*}$,  

{\em or\/} $p=q$;
\smallskip

\noindent $p\leq^{\alpha+1} q$\quad if and only if

{\em either\/} $p,q\in P^{\theta,\lambda}_\alpha$ and $p\leq^\alpha q$, 

{\em or\/} $q=\langle\zeta^*,\tau^*,n^*,u^*,\langle q_\xi,v_\xi:\xi<\theta
\rangle\rangle\in P^*_{\alpha+1}$, $p\in P^{\theta,\lambda}_\alpha$, and
$p\leq^\alpha q_\xi$ for some $\xi<\theta$,  

{\em or\/} $p=\langle\zeta^{**},\tau^*,n^*,u^*,\langle p_\xi,v_\xi:\xi<
\theta\rangle\rangle$, $q=\langle\zeta^*,\tau^*,n^*,u^*,\langle q_\xi,v_\xi:
\xi<\theta\rangle\rangle$ are from $P^*_{\alpha+1}$ and 
\[(\forall\xi<\theta)(p_\xi\leq^\alpha q_\xi\ \&\ u^{p_\xi}=u^{q_\xi}).\] 

\noindent [Again, it is easy to show that clauses (i)$_{\alpha+1}$--(iv)$_{
\alpha+1}$ are satisfied.]  
\medskip

After the construction is carried out we let 
\[\ptl=\bigcup\limits_{\alpha<\lambda} P^{\theta,\lambda}_\alpha\quad\mbox{
and }\quad{}\leq_{\rm pr}{}={}\bigcup\limits_{\alpha<\lambda}{\leq^\alpha_{
\rm pr}}\quad\mbox{ and }\quad{}\leq{}={}\bigcup\limits_{\alpha<\lambda}
{\leq^\alpha}.\]
One easily checks that $\leq_{\rm pr}$ is a partial order on $\ptl$ and that
the relation $\leq$ is transitive and reflexive, and that ${\leq_{\rm pr}} 
\subseteq {\leq}$.  

\begin{lemma}
\label{3.1}
Let $p,q\in\ptl$.
\begin{enumerate}
\item If $p\leq q$ then $\rht(p)\leq\rht(q)$, $u^p\subseteq u^q$ and $F^p=\{
f\rest u^p: f\in F^q\}$ (so $\bB^p$ is a subalgebra of $\bB^q$). If $p\leq q$
and $\rht(p)=\rht(q)$, then $q\leq p$.
\item For each $j\in u^p$, the set $\{\beta<\rht(p): h^p(j,\beta)<\theta\}$
is finite. 
\item If $p\leq_{\rm pr} q$ and $i\in u^p$, then $h^q(i,\beta)\geq\theta$
for all $\beta$ such that $\rht(p)\leq\beta<\rht(q)$.
\item If $i,j\in u^p$ are distinct, then there is $\beta<\rht(p)$ such that
$\theta\neq h^p(i,\beta)\neq h^p(j,\beta)\neq\theta$. 
\item For each finite set $X\subseteq\rht(p)$ there is $i\in u^p$ such that
\[\{\beta<\rht(p):h^p(i,\beta)<\theta\}=X.\] 
\item If $p\leq_{\rm pr} q$ then there is a $\leq_{\rm pr}$--increasing
sequence $\langle p_\xi:\xi\leq\rht(p)\rangle\subseteq\ptl$ such that
$p_{\rht(p)}=p$, $p_{\rht(q)}=q$ and $\rht(p_\xi)=\xi$ (for $\xi\leq\rht(p)$). 
(In particular, if $p\leq_{\rm pr} q$ and $\rht(p)=\rht(q)$ then $p=q$.)
\item If $\rht(p)=\gamma$ is a limit ordinal, $p=\langle p_\xi:\xi<\gamma
\rangle$, then for each $i\in u^p$ and $\xi<\gamma$:
\[i\in u^{p_\xi}\quad\mbox{ if and only if }\quad (\forall\zeta<\gamma)(\xi
\leq\zeta\ \Rightarrow\ h^p(i,\zeta)\geq\theta).\] 
\end{enumerate}
\end{lemma}

\begin{proof} 1)\quad Should be clear (an easy induction).
\smallskip 

\noindent 2)\quad Suppose that $p\in\ptl$ and $j\in u^p$ are a counterexample
with the minimal possible value of $\rht(p)$. Necessarily $\rht(p)$ is a limit
ordinal, $p=\langle p_\xi:\xi<\rht(p)\rangle$, $\rht(p_\xi)=\xi$ and $\zeta<
\xi<\rht(p)\ \Rightarrow\ p_\zeta\leq_{\rm pr} p_\xi$. Let $\xi<\rht(p)$ be
the first ordinal such that $j\in u^{p_\xi}$. By the choice of $p$, the set
$\{\beta\leq\xi: h^p(j,\beta)<\theta\}$ is finite, but clearly $h^p(j,\beta)
\geq\theta$ for all $\beta\in (\xi,\rht(p))$. 
\smallskip

\noindent 3)\quad An easy induction on $\rht(q)$ (with fixed $p$).
\smallskip

\noindent 4)\quad We show this by induction on $\rht(p)$. Suppose that $\rht(
p)=\alpha+1$, so $p=\langle\zeta^*,\tau^*,n^*,u^*,\langle p_\xi,v_\xi:\xi<
\theta\rangle\rangle$, and $i,j\in u^p$ are distinct. If $i,j\in u^{p_\xi}$
for some $\xi<\theta$, then by the inductive hypothesis we find $\beta<
\alpha$ such that 
\[\theta\neq h^p(i,\beta)=h^{p_\xi}(i,\beta)\neq h^{p_\xi}(j,\beta)=
h^p(j,\beta)\neq\theta.\]
If $i\in u^{p_\xi}\setminus u^*$, $j\in u^{p_\zeta}\setminus u^*$ and $\xi,
\zeta<\theta$ are distinct, then look at the definition of $h^p(i,\alpha)$,
$h^p(j,\alpha)$ -- these two values cannot be equal (and both are distinct
from $\theta$). Finally suppose that $\rht(p)$ is limit, so $p=\langle
p_\xi:\xi<\rht(p)\rangle$. Take $\xi<\rht(p)$ such that $i,j\in u^{p_\xi}$
and apply the inductive hypothesis to $p_\xi$ getting $\beta<\xi$ such that
$h^p(i,\beta)\neq h^p(j,\beta)$ (and both are not $\theta$).   
\smallskip

\noindent 5)\quad Again, it goes by induction on $\rht(p)$. First consider
a limit stage, and suppose that $\rht(p)=\gamma$ is a limit ordinal,
$X\in[\gamma]^{\textstyle{<}\omega}$ and $p=\langle p_\xi:\xi<\gamma
\rangle$. Let $\xi<\gamma$ be such that $X\subseteq\xi$. By the inductive
hypothesis we find $i\in u^{p_\xi}$ such that $\{\beta<\xi:h^p(i,\beta)<
\theta\}=X$. Applying clause (3) we may conclude that this $i$ is as
required. Now consider a successor case $\rht(p)=\alpha+1$. Let $p=\langle
\zeta^*,\tau^*,n^*,u^*,\langle p_\xi,v_\xi:\xi<\theta\rangle\rangle$, and
let $\xi<\theta$ be $\zeta^*$ if $\alpha\in X$, and be $\zeta^*+1$
otherwise. Apply the inductive hypothesis to $p_\xi$ and $X\cap \alpha$ to
get suitable $i\in u^{p_\xi}$, and note that this $i$ works for $p$ and $X$
too. 
\smallskip
  
\noindent 6), 7)\quad Straightforward. 
\end{proof}

\begin{definition}
\label{defiso}
We say that conditions $p,q\in\ptl$ are {\em isomorphic\/} if $\rht(p)=
\rht(q)$, $\otp(u^p)=\otp(u^q)$, and if $H:u^p\longrightarrow u^q$ is the
order isomorphism, then for every $\beta<\rht(p)$
\[(\forall j\in u^p)(h^p(j,\beta)=h^q(H(j),\beta)\ \&\ g^p(j,\beta)=g^p(H(j),
\beta)).\]
[In this situation we may say that $H$ is the isomorphism from $p$ to $q$.]
\end{definition}

\begin{lemma}
\label{3x1}
Suppose that $q_0,q_1\in\ptl$ are isomorphic conditions and $H$ is the
isomorphism from $q_0$ to $q_1$. 
\begin{enumerate}
\item If $\rht(q_0)=\rht(q_1)=\gamma$ is a limit ordinal, $q_\ell=\langle
q^\ell_\xi: \xi<\gamma\rangle$ (for $\ell<2$), then $H\rest u^{q_\xi^0}$ is
an isomorphism from $q^0_\xi$ to $q^1_\xi$. 
\item If $\rht(q_0)=\rht(q_1)=\alpha+1$, $\alpha<\lambda$, and $q_\ell=
\langle\zeta^*_\ell,\tau^*_\ell,n^*_\ell,u^*_\ell,\langle q^\ell_\xi,
v^\ell_\xi:\xi<\theta\rangle\rangle$ (for $\ell<2$), then $\zeta^*_0=
\zeta^*_1$, $\tau^*_0=\tau^*_1$, $n^*_0=n^*_1$, $H\rest u^{q_\xi^0}$ is an
isomorphism from $q^0_\xi$ to $q^1_\xi$ and $H[v^0_\xi]=v^1_\xi$ (for $\xi<
\theta$).   
\item $F^{q_0}=\{f\comp H:f\in F^{q_1}\}$.
\item Assume $p_0\leq q_0$. Then there is a unique condition $p_1\leq q_1$
such that $H\rest u^{p_0}$ is the isomorphism from $p_0$ to $p_1$.\\
\relax [The condition $p_1$ will be called $H(p_0)$.]
\end{enumerate}
\end{lemma}

\begin{proof}
1), 2)\quad Straightforward (for (1) use Lemma \ref{3.1}(7)).\\
3), 4)\quad Easy inductions on $\rht(q_0)$ using (1), (2) above.
\end{proof}

\begin{definition}
\label{defptran}
By induction on $\alpha<\lambda$, for conditions $p,q\in P^{\theta,
\lambda}_\alpha$ such that $p\leq^\alpha q$, we define {\em the
$p$--transformation $T_p(q)$ of $q$}. 
\begin{itemize}
\item If $\alpha=0$ (so necessarily $p=q$) then $T_p(q)=p$.
\item Assume that $\rht(q)=\alpha+1$, $q=\langle\zeta^*,\tau^*,n^*,u^*,
\langle q_\xi,v_\xi:\xi<\theta\rangle\rangle$. 

If $p\leq q_\xi$ for some $\xi<\theta$, then let $\xi^*$ be such that $p\leq
q_{\xi^*}$. Next for $\xi<\theta$ let $q_\xi'=T_{H_{\xi^*, \xi}(p)}(q_\xi)$,
where $H_{\xi^*,\xi}$ is the isomorphism from $q_{\xi^*}$ to $q_\xi$. Define
$T_p(q)=\langle\xi^*,\tau^*,n^*,u^*,\langle q_\xi',v_\xi: \xi<\theta\rangle
\rangle$. 
 
Suppose now that $p=\langle\zeta^{**},\tau^*,n^*,u^*,\langle p_\xi,v_\xi:\xi 
<\theta\rangle\rangle$ and $u^{p_\xi}=u^{q_\xi}$, $p_\xi\leq q_\xi$ (for
$\xi<\theta$). Let $q_\xi'=T_{p_\xi}(q_\xi)$ and put $T_p(q)=\langle
\zeta^{**},\tau^*,n^*,u^*,\langle q_\xi',v_\xi:\xi<\theta\rangle\rangle$.
\item Assume now that $\rht(q)$ is a limit ordinal and $q=\langle q_\xi:\xi<
\rht(q)\rangle$. 

If $\rht(p)<\rht(q)$ then $p\leq q_\vare$ for some $\vare<\rht(q)$, and we
may choose $q_\xi'$ (for $\xi<\rht(q)$) such that $\rht(q_\xi')=\xi$, $\xi<
\xi'<\rht(q)\ \Rightarrow\ q_\xi'\leq_{\rm pr} q_{\xi'}'$, and $q_\zeta'
=T_p(q_\zeta)$ for $\zeta\in [\vare,\rht(q))$. Next we let $T_p(q)=\langle
q_\zeta':\zeta<\theta\rangle$. 

If $\rht(p)=\rht(q)$, $p=\langle p_\xi:\xi<\rht(p)\rangle$ and $p_\xi\leq
q_\xi$ for $\xi>\delta$ (for some $\delta<\rht(p)$) then we define
$T_p(q)=p$. 
\end{itemize}
\end{definition}
To show that the definition of $T_p(q)$ is correct one proves inductively
(parallely to the definition of the $p$--transformation of $q$) the following
facts. 

\begin{lemma}
\label{3y1}
Assume $p,q\in\ptl$, $p\le q$. Then:
\begin{enumerate}
\item $T_p(q)\in\ptl$, $u^{T_p(q)}=u^q$, $\rht(T_p(q))=\rht(q)$,
\item $p\leq_{\rm pr} T_p(q)\leq q\leq T_p(q)$,
\item $\rht(p)=\rht(q)\ \Rightarrow\ T_p(q)=p$,
\item if $q'\in\ptl$ is isomorphic to $q$ and $H:u^q\longrightarrow u^{q'}$
is the isomorphism from $q$ to $q'$, then $H$ is the isomorphism from
$T_p(q)$ to $T_{H(p)}(q')$,
\item if $q\leq_{\rm pr}q'$ then $T_p(q)\leq_{\rm pr} T_p(q')$. 
\end{enumerate}
\end{lemma}

\begin{proposition}
\label{3y2}
Every $\leq_{\rm pr}$--increasing chain in $\ptl$ of length $<\lambda$ has
a $\leq_{\rm pr}$--upper bound, that is the partial order $(\ptl,\leq_{\rm
pr})$ is $(<\lambda)$--closed.
\end{proposition}

Let us recall that a forcing notion $(\bQ,\leq)$ is {\em
$({<}\lambda)$--strategically closed\/} if the second player has a winning
strategy in the following game $\Game_\lambda(\bQ)$.

The game $\Game_\lambda(\bQ)$ lasts $\lambda$ moves. The first player starts
with choosing a condition $p^*\in\bQ$. Later, in her $i^{\rm th}$ move, the
first player chooses an open dense subset $D_i$ of $\bQ$. The second player
(in his $i^{\rm th}$ move) picks a condition $p_i\in\bQ$ so that $p_0\geq
p^*$, $p_i\in D_i$ and $p_i\geq p_j$ for all $j<i$. The second player looses
the play if for some $i<\lambda$ he has no legal move. 

It should be clear that $({<}\lambda)$--strategically closed forcing notions
do not add sequences of ordinals of length less than $\lambda$. The reader
interested in this kind of properties of forcing notions and iterating them
is referred to \cite{Sh:587}, \cite{Sh:667}. 

\begin{proposition}
\label{3.3}
Assume that $\theta<\lambda$ are regular cardinals, $\lambda^{<\lambda}=
\lambda$. Then $(\ptl,\leq)$ is a $(<\lambda)$--strategically closed
$\lambda^+$--cc forcing notion.
\end{proposition}

\begin{proof}
It follows from Lemma \ref{3y1}(2) that if $D\subseteq\ptl$ is an open dense
set, $p\in\ptl$, then there is a condition $q\in D$ such that $p\leq_{\rm
pr} q$. Therefore, to win the game $\Game_\lambda(\ptl)$, the second player
can play so that the conditions $p_i$ that he chooses are $\leq_{\rm 
pr}$--increasing, and thus there are no problems with finding $\leq_{\rm
pr}$--bounds (remember Proposition \ref{3y2}).

Now, to show that $\ptl$ is $\lambda^+$--cc, suppose that $\langle p_\delta:
\delta<\lambda^+\rangle$ is a sequence of distinct conditions from
$\ptl$. We may find a set $A\in [\lambda^+]^{\textstyle \lambda^+}$ such
that 
\begin{itemize}
\item conditions $\{p_\delta:\delta\in A\}$ are pairwise isomorphic,
\item the family $\{u^{p_\delta}:\delta\in A\}$ forms a $\Delta$--system
with heart $u^*$,
\item if $\delta_0<\delta_1$ are from $A$ then 
\[\sup(u^*)<\min(u^{p_{\delta_0}}\setminus u^*)\leq \sup(u^{p_{\delta_0}}
\setminus u^*)<\min(u^{p_{\delta_0}}\setminus u^*).\]
\end{itemize}
Take an increasing sequence $\langle\delta_\xi:\xi<\theta\rangle$ of
elements of $A$, let $\tau^*={\bf 1}$, $v_\xi=\emptyset$ (for $\xi<\theta$),
and look at $p=\langle 0,\tau^*,0,u^*,\langle p_{\delta_\xi},v_\xi:\xi<
\theta\rangle\rangle$. It is a condition in $\ptl$ stronger than all
$p_{\delta_\xi}$'s. 
\end{proof}

\begin{definition}
\label{defcompo}
By induction on $\rht(p)$ we define {\em $\alpha$--components of $p$\/} (for
$p\in\ptl$, $\alpha\leq\rht(p)$). 
\begin{itemize}
\item First we declare that the only $\rht(p)$--component of $p$ is the $p$
itself.  
\item If $\rht(p)=\beta+1$, $p=\langle\zeta^*,\tau^*,n^*,u^*,\langle p_\xi,
v_\xi:\xi<\theta\rangle\rangle$ and $\alpha=\beta$, then $\alpha$--components
of $p$ are $p_\xi$ (for $\xi<\theta$); if $\alpha<\beta$, then
$\alpha$--components of $p$ are those $q$ which are $\alpha$--components of
$p_\xi$ for some $\xi<\theta$. 
\item If $\rht(p)$ is a limit ordinal, $p=\langle p_\xi:\xi<\rht(p)\rangle$
and $\alpha<\rht(p)$, then $\alpha$--components of $p$ are
$\alpha$--components of $p_\xi$ for $\xi\in [\alpha,\rht(p))$.
\end{itemize}
\end{definition}

\begin{lemma}
\label{3.4x}
Assume $p\in \ptl$ and $\alpha<\rht(p)$. 
\begin{enumerate}
\item If $q$ is an $\alpha$--component of $p$ then $q\leq p$, $\rht(q)=
\alpha$, and for all $j_0,j_1\in u^q$ and every $\beta\in [\alpha,\rht(p))$: 
\[h^p(j_0,\beta)\neq\theta\ \&\ h^p(j_1,\beta)\neq \theta\quad\Rightarrow
\quad h^p(j_0,\beta)=h^p(j_1,\beta).\]
Moreover, for each $i\in u^p$ there is a unique $\alpha$--component $q$ of
$p$ such that $i\in u^q$ and
\[(\forall j\in u^q)(\forall\beta\in [\alpha,\rht(p)))(h^p(i,\beta)\geq
\theta\ \Rightarrow\ h^p(j,\beta)\geq\theta).\]
\item If $H$ is an isomorphism from $p$ onto $p'\in\ptl$, and $q$ is an
$\alpha$--component of $p$, then $H(q)$ is an $\alpha$--component of
$p'$. If $q_0,q_1$ are $\alpha$--components of $p$ then $q_0,q_1$ are
isomorphic.    
\item There is a unique $\alpha$--component $q$ of $p$ such that $q\leq_{\rm
pr} p$.
\end{enumerate}
\end{lemma}

\begin{proof}
Easy inductions on $\rht(p)$.
\end{proof}

\begin{definition}
\label{closed}
By induction on $\rht(p)$ we define when a set $Z\subseteq \lambda$ is
$p$--closed for a condition $p\in\ptl$.  
\begin{itemize}
\item If $\rht(p)=0$ then every $Z\subseteq\lambda$ is $p$--closed; 
\item if $\rht(p)$ is limit, $p=\langle p_\xi:\xi<\rht(p)\rangle$, then
$Z$ is $p$--closed provided it is $p_\xi$--closed for each $\xi<\rht(p)$;
\item if $\rht(p)=\alpha+1$, $p=\langle\zeta^*,\tau^*,n^*,u^*,\langle
p_\xi,v_\xi:\xi<\theta\rangle\rangle$ and $\alpha\notin Z$, then $Z$ is
$p$--closed whenever it is $p_{\zeta^*}$--closed;
\item if $\rht(p)=\alpha+1$, $p=\langle\zeta^*,\tau^*,n^*,u^*,\langle
p_\xi,v_\xi:\xi<\theta\rangle\rangle$ and $\alpha\in Z$, then $Z$ is
$p$--closed provided it is $p_{\zeta^*}$--closed and 
\[\{\beta<\alpha:(\exists j\in v_{\zeta^*}\cup\{\min(u^{p_{\zeta^*}}
\setminus u^*)\})(h^{p_{\zeta^*}}(j,\beta)<\theta)\}\subseteq Z.\]
\end{itemize}
\end{definition}

\begin{lemma}
\label{3.5.1}
\begin{enumerate}
\item If $p\in\ptl$ and $w\in [\rht(p)]^{\textstyle<\omega}$, then there is
a finite $p$--closed set $Z\subseteq\rht(p)$ such that $w\subseteq Z$.
\item If $p,q\in\ptl$ are isomorphic and $Z$ is $p$--closed, then $Z$ is
$q$--closed. If $Z$ is $p$--closed, $\alpha<\rht(p)$ and $p^*$ is an
$\alpha$--component of $p$, then $Z\cap\alpha$ is $p^*$--closed. 
\end{enumerate}
\end{lemma}

\begin{proof}
Easy inductions on $\rht(p)$ (remember Lemma \ref{3.1}(2)).  
\end{proof}

\begin{definition}
\label{UUpsilon}
Suppose that $p\in\ptl$ and $Z\subseteq \rht(p)$ is a finite $p$--closed
set. Let $Z=\{\alpha_0,\ldots,\alpha_{k-1}\}$ be the increasing
enumeration. 
\begin{enumerate}
\item We define 
\[U[p,Z]\stackrel{\rm def}{=}\{j\in u^p: (\forall\beta<\rht(p))(h^p(j,\beta) 
<\theta\ \Rightarrow\ \beta\in Z)\}.\]
\item We let
\[\Upsilon_p(Z)=\langle\zeta_\ell,\tau_\ell,n_\ell,\langle g_\ell,h^\ell_0,
\ldots,h^\ell_{n_\ell-1}\rangle:\ell<k\rangle,\]
where, for $\ell<k$, $\zeta_\ell$ is an ordinal below $\theta$, $\tau_\ell$
is a Boolean term, $n_\ell<\omega$ and $g_\ell,h^\ell_0,\ldots,h^\ell_{
n_\ell-1}:\ell\longrightarrow 2$, and they all are such that for every
(equivalently: some) $\alpha_\ell+1$--component $q=\langle\zeta^*,\tau^*,
n^*,u^*,\langle q_\xi,v_\xi:\xi<\theta\rangle\rangle$ of $p$ we have:\quad
$\zeta_\ell=\zeta^*$, $\tau_\ell=\tau^*$, $n_\ell=n^*$ and if $v_\xi=\{j_0,
\ldots,j_{n_\ell-1}\}$ (the increasing enumeration) then   
\[(\forall m<n_\ell)(\forall\ell'<\ell)(h^\ell_m(\ell')=h^q(j_m,
\alpha_{\ell'})),\]
and if $i_0=\min(u^{q_{\zeta^*}}\setminus u^*)$ then $(\forall\ell'<\ell)(
g_\ell(\ell')=h^q(i_0,\alpha_{\ell'}))$. (Note that $\zeta_\ell,\tau_\ell,
n_\ell$, $g_\ell,h^\ell_0,\ldots,h^\ell_{n_\ell-1}$ are well-defined by
Lemma \ref{3.4x}. Necessarily, for all $m<n_\ell$ and $\beta\in
\alpha_\ell\setminus Z$ we have $h^q(i_0,\beta),h^q(j_m,\beta)\geq\theta$;
remember that $Z$ is $p$--closed.)  
\end{enumerate}
\end{definition}

Note that if $Z\subseteq\rht(p)$ is a finite $p$--closed set,
$\alpha=\max(Z)$ and $p^*$ is the $\alpha+1$--component of $p$ satisfying
$p^*\leq_{\rm pr} p$ (see \ref{3.4x}(3)), then $U[p,Z]\subseteq u^{p^*}$. 

\begin{lemma}
\label{3.6.x}
Suppose that $p\in\ptl$ and $Z_0,Z_1\subseteq\rht(p)$ are finite $p$--closed 
sets such that $\Upsilon_p(Z_0)=\Upsilon_p(Z_1)$. Then $\otp(U[p,Z_0])=
\otp(U[p,Z_1])$, and the order preserving isomorphism
$\pi:U[p,Z_0]\longrightarrow U[p,Z_1]$ satisfies
\begin{enumerate}
\item[$(\otimes)$] \quad $(\forall\ell<k)(h^p(i,\alpha^0_\ell)=h^p(\pi(i),
\alpha^1_\ell))$, 

where $\{\alpha^x_0,\ldots,\alpha^x_{k-1}\}$ is the increasing enumeration
of $Z_x$ (for $x=0,1$). 
\end{enumerate}
\end{lemma}

\begin{proof}
We prove this by induction on $|Z_0|=|Z_1|$ (for all $p,Z_0,Z_1$ satisfying
the assumptions).
\medskip

\noindent {\sc Step} $|Z_0|=|Z_1|=1$; $Z_0=\{\alpha^0_0\}$,
$Z_1=\{\alpha^1_0\}$.\\
Take the $\alpha^x_0+1$--component $q_x$ of $p$ such that $q_x\leq_{\rm pr}
p$. Then, for $x=0,1$, $q_x=\langle\zeta,\tau,n, u^x,\langle q^x_\xi,
v^x_\xi:\xi<\theta\rangle\rangle$, and for each $i\in v^x_\xi$, $\beta<
\alpha^x_0$ we have $h^{q^x_\xi}(i,\beta)\geq\theta$. Also, if $i^x_0=
\min(u^{q^x_\zeta}\setminus u^x)$ and $\beta<\alpha^x_0$, then
$h^{q^x_\zeta}(i^x_0,\beta)\geq\theta$. Consequently, $n=|v^x_\xi|\leq 1$,
and if $n=1$ then $\{i^x_0\}=v^x_\zeta$ (remember Lemma \ref{3.1}(4)). 
Moreover, 
\[U[p,Z_x]=U[q_x,Z_x]=\{H^x_{\xi,\zeta}(i^x_0):\xi<\theta\},\]
where $H^x_{\xi,\zeta}$ is the isomorphism from $q^x_\zeta$ to
$q^x_\xi$. Now it should be clear that the mapping $\pi:H^0_{\xi,\zeta}(
i^0_0)\mapsto H^1_{\xi,\zeta}(i^1_0):U[p,Z_0]\longrightarrow U[p,Z_1]$ is
the order preserving isomorphism  (remember clause $(\gamma)$ of the
definition of $P^*_{\alpha+1}$), and it has the property described in
$(\otimes)$. 
\medskip

\noindent {\sc Step} $|Z_0|=|Z_1|=k+1$; $Z_0=\{\alpha^0_0,\ldots,\alpha^0_k
\}$, $Z_1=\{\alpha^1_0,\ldots,\alpha^1_k\}$.\\
Let 
\[\Upsilon_p(Z_0)=\Upsilon_p(Z_1)=\langle\zeta_\ell,\tau_\ell,n_\ell,\langle
g_\ell,h^\ell_0,\ldots,h^\ell_{n_\ell-1}\rangle:\ell\leq k\rangle.\]
For $x=0,1$, let $q_x=\langle\zeta,\tau,n, u^x,\langle q^x_\xi,v^x_\xi:\xi<
\theta\rangle\rangle$ be the $\alpha^x_k+1$--component of $p$ such that $q_x
\leq_{\rm pr} p$. The sets $Z_x\cap\alpha^x_k$ (for $x=0,1$) are
$q^x_\xi$--closed for every $\xi<\theta$, and clearly $\Upsilon_p(Z_0\cap
\alpha^0_k)= \Upsilon_p(Z_1\cap\alpha^1_k)$. Hence, by the inductive
hypothesis, $\otp(U[q^0_\xi,Z_0\setminus\{\alpha^0_k\}])=\otp(U[q^1_\xi,Z_1
\setminus\{\alpha^1_k\}])$ (for each $\xi<\theta$), and the order preserving
mappings $\pi_\xi:U[q^0_\xi,Z_0\setminus\{ \alpha^0_k\}]\longrightarrow
U[q^1_\xi, Z_1\setminus\{\alpha^1_k\}]$ satisfy the demand in $(\otimes)$.
Let $i^x_\xi=\min(u^{q^x_\xi}\setminus u^x)$. Then, as $q^x_\xi$ and
$q^x_\zeta$ are isomorphic and the isomorphism is the identity on $u^x$, we
have $(\forall\ell<k)(h^p(i^x_\xi,\alpha_\ell^x)=g_k(\ell))$. Hence $\pi_\xi
(i^0_\xi)=i^1_\xi$, and therefore $\pi_\xi[u^0\cap U[q^0_\xi,Z_0\setminus
\{\alpha^0_k\}]]=u^1\cap U[q^1_\xi,Z_1\setminus\{\alpha^1_k\}]$. But since
the mappings $\pi_\xi$ are order preserving, the last equality implies that
$\pi_\xi\restriction (u^0\cap U[q^0_\xi,Z_0\setminus\{\alpha^0_k\}])=
\pi_\zeta\restriction (u^0\cap U[q^0_\zeta,Z_0\setminus\{\alpha^0_k\}])$,
and hence $\pi=\bigcup\limits_{\xi<\theta}\pi_\xi$ is a function, and it is
an order isomorphism from $U[q_0,Z_0]=U[p,Z_0]$ onto $U[q_1,Z_1]=U[p,Z_1]$
satisfying $(\otimes)$.  
\end{proof}

\section{The algebra and why it is OK (in $\bV^{\ptl}$)}
Let $\dbtl$ and $\dot{U}$ be $\ptl$--names such that 
\[\forces_{\ptl}\mbox{`` }\dbtl=\bigcup\{\bB^p:p\in \Gamma_{\ptl}\}\mbox{
''}\quad\mbox{ and }\quad\forces_{\ptl}\mbox{`` }\dot{U}=\bigcup\{u^p:p\in 
\Gamma_{\ptl}\}\mbox{ ''.}\] 
Note that $\dot{U}$ is (a name for) a subset of $\lambda^+$. Let $\dot{F}$
be a $\ptl$--name such that 
\[\forces_{\ptl}\mbox{`` }\dot{F}=\{f\in 2^{\textstyle\dot{U}}: (\forall
p\in\Gamma_{\ptl})(f\restriction u^p\in\dot{F}^p)\}\mbox{ ''.}\]

\begin{proposition}
\label{3.4}
Assume $\theta<\lambda$ are regular, $\lambda^{<\lambda}=\lambda$. Then in
$\bV^{\ptl}$: 
\begin{enumerate}
\item $\dot{F}$ is a non-empty closed subset of $2^{\textstyle\dot{U}}$, and
$\dbtl$ is the Boolean algebra generated $\bB_{(\dot{U},\dot{F})}$ (see
Definition \ref{0.C});
\item $|\dot{U}|=|\dbtl|=\lambda^+$;
\item For every subalgebra $\bB\subseteq\dbtl$ of size $\lambda^+$ we have
$\Dep^+(\bB)>\theta$. 
\end{enumerate}
\end{proposition}

\begin{proof} 
2)\quad Note that if $p\in\ptl$, $\sup(u^p)<j<\lambda^+$ then there is a
condition $q\geq p$ such that $j\in u^q$. Hence $\forces|\dot{U}|=
\lambda^+$. To show that, in $\bV^{\ptl}$, the algebra $\dbtl$ is of size
$\lambda^+$ it is enough to prove the following claim. 
\begin{claim}
\label{3.4.1}
Let $p\in\ptl$, $j\in u^p$. Then $x_j\notin\langle x_i:i\in j\cap u^p
\rangle_{\bB^p}$.
\end{claim}

\begin{proof}[Proof of the claim]
Suppose not, and let $p,j$ be a counterexample with the smallest possible
$\rht(p)$. Necessarily, $\rht(p)$ is a successor ordinal, say
$\rht(p)=\alpha+1$. So let $p=\langle\zeta^*,\tau^*, n^*,u^*,\langle
p_\xi,v_\xi:\xi<\theta\rangle\rangle$ and suppose that $v\in [u^p\cap
j]^{\textstyle <\omega}$ is such that $x_j\in \langle x_i: i\in
v\rangle_{\bB^p}$. If $j\in u^*$ then $v\subseteq u^*$ and we immediately
get a contradiction (applying the inductive hypothesis to $p_{\zeta^*}$). So
let $\xi<\theta$ be such that $j\in u^{p_\xi}\setminus u^*$. We know that
$x_j\notin\langle x_i:i\in u^*\cup (v\cap u^{p_\xi})\rangle_{\bB^{p_\xi}}$
(remember clause $(\gamma)$ of the definition of $P^*_{\alpha+1}$), so we may
take functions $f_0,f_1\in F^{p_\xi}$ such that $f_0\rest (u^*\cup (v\cap
u^{p_\xi}))=f_1\rest (u^*\cup (v\cap u^{p_\xi}))$, $f_0(j)=0$, $f_1(j)=1$. 
Let $g_0,g_1:u^p\longrightarrow 2$ be such that $g_\ell\rest u^{p_\xi}=
f_\ell$, $g_\ell\rest u^{p_\zeta}=f_0\comp H_{\zeta,\xi}$ for $\zeta\neq\xi$
(where $H_{\zeta,\xi}$ is the order isomorphism from $u^{p_\zeta}$ to
$u^{p_\xi}$). Now one easily checks that $g_0,g_1\in F^p$ (remember the
definition of the term $\sigma_{\rm maj}$). By our choices, $g_0(i)=g_1(i)$
for all $i\in v$, and $g_0(j)\neq g_1(j)$, and this is a clear contradiction
with the choice of $i$ and $v$.
\end{proof}

\noindent 3)\quad Suppose that $\langle\dot{a}_\xi:\xi<\lambda^+\rangle$ is
a $\ptl$--name for a $\lambda^+$--sequence of distinct members of $\dbtl$
and let $p\in\ptl$. Applying standard cleaning procedures we find a set
$A\subseteq\lambda^+$ of the order type $\theta$, an ordinal $\alpha<
\lambda$ and $\tau^*,n^*,u^*$ and $\langle p_\xi,v_\xi:\xi\in A\rangle$ such
that $p\leq p_\xi$, $\rht(p_\xi)=\alpha$, $p_\xi\forces\dot{a}_\xi=
\tau^*(x_i:i\in v_\xi)$ and 
\[q\stackrel{\rm def}{=}\langle 0,\tau^*,n^*,u^*,\langle p_\xi,v_\xi:\xi\in
A\rangle\rangle\in P^*_{\alpha+1},\]
where $A$ is identified with $\theta$ by the increasing enumeration (so we
will think $A=\theta$). For $\xi<\theta$ let $\tau_\xi=\tau^*(x_i: i\in
v_\xi)\in\bB^{p_\xi}$. Since $\dot{a}_\xi$ were (forced to be) distinct we
know that $\bB^q\models\tau_\xi\neq\tau_\zeta$ for distinct $\xi,\zeta$. 
Hence $\tau_\xi\notin\langle x_i:i\in u^*\rangle_{\bB^{p_\xi}}$ (for each
$\xi$) and therefore we may find functions $f^0_\xi,f^1_\xi\in F^{p_\xi}$
such that $f^0_\xi\rest u^*=f^1_\xi\rest u^*$, and $f^0_\xi(\tau_\xi)=0$,
$f^1_\xi(\tau_\xi)=1$, and if $\xi<\zeta<\theta$, and $H_{\xi,\zeta}$ is the
isomorphism from $p_\xi$ to $p_\zeta$, then $f^\ell_\xi=f^\ell_\zeta\comp
H_{\xi,\zeta}$. Now fix $\xi<\zeta<\theta$ and let  
\[g\stackrel{\rm def}{=}\bigcup_{\alpha\leq 3\cdot\xi+2} f^0_\alpha\cup
\bigcup_{3\cdot\xi+2<\alpha<\theta} f^1_\alpha.\]
It should be clear that $g$ is a function from $u^q$ to $2$, and moreover
$g\in F^q$. Also easily
\[g(\sigma_{\rm maj}(\tau_{3\cdot\xi},\tau_{3\cdot\xi+1},\tau_{3\cdot\xi+ 
2}))=0\ \mbox{ and }\ g(\sigma_{\rm maj}(\tau_{3\cdot\zeta},\tau_{3\cdot
\zeta+1},\tau_{3\cdot\zeta+2}))\}=1.\]
Hence we may conclude that
\[\bB^q\models\sigma_{\rm maj}(\tau_{3\cdot\xi},\tau_{3\cdot\xi+1},\tau_{3
\cdot\xi+2})<\sigma_{\rm maj}(\tau_{3\cdot\zeta},\tau_{3\cdot\zeta+1},\tau_{3
\cdot\zeta+2})\]
for $\xi<\zeta<\theta$ (remember the definition of $F^q$ and Proposition
\ref{0.D}). Consequently we get $q\forces\Dep^+(\langle\dot{a}_\xi:\xi<
\lambda^+\rangle_{\dbtl})>\theta$, finishing the proof. 
\end{proof}

\begin{theorem}
\label{3.5}
Assume $\theta<\lambda$ are regular, $\lambda=\lambda^{<\lambda}$. Then
$\forces_{\ptl}\Dep(\dbtl)=\theta$. 
\end{theorem}

\begin{proof}
By Proposition \ref{3.4} we know that $\forces\Dep^+(\dbtl)>\theta$, so what
we have to show is that there are no increasing sequences of length
$\theta^+$ of elements of $\dbtl$. We will show this under an additional
assumption that $\theta^+ <\lambda$ (after the proof is carried out, it will 
be clear how one modifies it to deal with the case $\lambda=\theta^+$). Due
to this additional assumption, and since the forcing notion $\ptl$ is
$(<\lambda)$--strategically closed (by Proposition \ref{3.3}), it is enough
to show that $\Dep(\bB^p)\leq\theta$ for each $p\in\ptl$. 

So suppose that $p\in\ptl$ is such that $\Dep(\bB^p)\geq\theta^+$. Then we
find a Boolean term $\tau$, an integer $n$ and sets $w_\rho\in [u^p]^{
\textstyle n}$ (for $\rho<\theta^+$) such that
\[\rho_0<\rho_1<\theta^+\quad\Rightarrow\quad\bB^p\models\tau(x_i:i\in w_{
\rho_0})<\tau(x_i:i\in w_{\rho_1}).\]
For each $\rho<\theta^+$ use Lemma \ref{3.5.1} to choose a finite
$p$--closed set $Z_\rho\subseteq\rht(p)$ containing the set 
\[\{\beta<\rht(p): (\exists j\in w_\rho)(h^p(j,\beta)<\theta)\}.\] 
Look at $\Upsilon_p(Z_\rho)$ (see Definition \ref{UUpsilon}). There are only
$\theta$ possibilities for the values of $\Upsilon_p(Z_\rho)$, so we find
$\rho_0<\rho_1<\theta^+$ such that   
\begin{enumerate}
\item[(i)] $|Z_{\rho_0}|=|Z_{\rho_1}|$, $\Upsilon_p(Z_{\rho_0})=\Upsilon_p( 
Z_{\rho_1})=\langle\zeta_\ell,\tau_\ell,n_\ell,\langle g_\ell,h^\ell_0,
\ldots,h^\ell_{n_\ell-1}\rangle:\ell<k\rangle$,
\item[(ii)] if $\pi^*:Z_{\rho_0}\longrightarrow Z_{\rho_1}$ is the order
isomorphism then $\pi^*\rest Z_{\rho_0}\cap Z_{\rho_1}$ is the identity on
$Z_{\rho_0}\cap Z_{\rho_1}$,
\item[(iii)] if $\pi:U[p,Z_{\rho_0}]\longrightarrow U[p,Z_{\rho_1}]$ is the
order isomorphism, then $\pi[w_{\rho_0}]=w_{\rho_1}$.
\end{enumerate}
Note that, by Lemma \ref{3.6.x}, $\otp(U[p,Z_{\rho_0}])=\otp(U[p,
Z_{\rho_1}])$ and the order isomorphism $\pi$ satisfies 
\[(\forall j\in U[p,Z_{\rho_0}])(\forall\beta\in Z_{\rho_0})(h^p(j,\beta)=
h^p(\pi(j),\pi^*(\beta))),\]
and hence $\pi$ is the identity on $U[p,Z_{\rho_0}]\cap U[p,Z_{\rho_1}]$
(remember Lemma \ref{3.1}).  

For a function $f\in F^p$ let $G^{\rho_0}_{\rho_1}(f):u^p\longrightarrow 2$
be defined by 
\[G^{\rho_0}_{\rho_1}(f)(j)=\left\{\begin{array}{lll}
f(\pi(j))&\mbox{ if }& j\in U[p,Z_{\rho_0}],\\
f(\pi^{-1}(j))&\mbox{ if }& j\in U[p,Z_{\rho_1}]\setminus U[p,\rho_0],\\
0& &\mbox{otherwise}.
				   \end{array}\right.\]
\begin{claim}
\label{3.5.2}
For each $f\in F^p$, $G^{\rho_0}_{\rho_1}(f)\in F^p$.
\end{claim}

\begin{proof}[Proof of the claim]
By induction on $\alpha\leq\rht(p)$ we show that for each
$\alpha$--component $q$ of $p$, the restriction $G^{\rho_0}_{\rho_1}(f)\rest 
u^q$ is in $F^q$.  

If $\alpha$ is limit, we may easily use the inductive hypothesis to show
that, for any $\alpha$--component $q$ of $p$, $G^{\rho_0}_{\rho_1}(f)\rest
u^q\in F^q$.   

Assume $\alpha=\beta+1$ and let $q=\langle\zeta^*,\tau^*,n^*,u^*,\langle
q_\xi,v_\xi:\xi<\theta\rangle\rangle$ be an $\alpha$--component of $p$. We
will consider four cases.
\medskip

\noindent{\em Case 1:}\quad $\beta\notin Z_{\rho_0}\cup Z_{\rho_1}$.\\
Then $(U[p,Z_{\rho_0}]\cup U[p,Z_{\rho_1}])\cap u^q\subseteq
u^{q_{\zeta^*}}$ and $G^{\rho_0}_{\rho_1}(f)\rest (u^{q_\xi}\setminus u^*)
\equiv 0$ for each $\xi\neq\zeta^*$. Since, by the inductive hypothesis,
$G^{\rho_0}_{\rho_1}(f)\rest u^{q_\xi}\in F^{q_\xi}$ for each $\xi<\theta$,
we may use the definition of $P^*_{\beta+1}$ and conclude that
$G^{\rho_0}_{\rho_1}(f)\rest u^q\in F^q$ (remember the definition of the 
term $\sigma_{\rm maj}$). 
\medskip

\noindent{\em Case 2:}\quad $\beta\in Z_{\rho_0}\setminus Z_{\rho_1}$.\\
Let $Z_{\rho_0}=\{\alpha_0,\ldots,\alpha_{k-1}\}$ be the increasing
enumeration. Then $\beta=\alpha_\ell$ for some $\ell<k$ and $\zeta^*=
\zeta_\ell$, $\tau^*=\tau_\ell$, $n^*=n_\ell$. Moreover, if $v_\xi= 
\{j^\xi_0,\ldots,j^\xi_{n_\ell-1}\}$ (the increasing enumeration),
$\xi<\theta$, then for $m<n_\ell$: 
\[(\forall\ell'<\ell)(h^\ell_m(\alpha_{\ell'})=h^q(j^\xi_m,\alpha_{\ell'}))
\quad\mbox{and}\quad(\forall\gamma\in\beta\setminus Z_{\rho_0})(h^q(j^\xi_m, 
\gamma)\geq\theta).\] 
Note that $U[p,Z_{\rho_1}]\cap u^q\subseteq u^{q_{\zeta^*}}$, so if
$U[p,Z_{\rho_0}]\cap u^q=\emptyset$, then we may proceed as in the previous
case. Therefore we may assume that $U[p,Z_{\rho_0}]\cap u^q\neq\emptyset$. 
So, for each $\gamma\in Z_{\rho_0}\setminus\alpha$ we may choose $i_\gamma
\in U[p,Z_{\rho_0}]\cap u^q$ such that 
\[(\forall i\in U[p,Z_{\rho_0}]\cap u^q)(h^p(i,\gamma)\neq\theta\
\Rightarrow\ h^p(i,\gamma)=h^p(i_\gamma,\gamma))\]
(remember Lemma \ref{3.4x}(1)).  Let $i^*=\max\{i_\gamma:\gamma\in
Z_{\rho_0}\setminus\alpha\}$ (if $\beta=\max(Z_{\rho_0})$, then let $i^*$ be 
any element of $U[p,Z_{\rho_0}]\cap u^q$). Note that then
\[(\forall i\in U[p,Z_{\rho_0}]\cap u^q)(\forall\gamma\in Z_{\rho_0}
\setminus\alpha)(h^p(i,\gamma)\neq\theta\ \Rightarrow\ h^p(i,\gamma)=h^p(
i^*,\gamma))\] 
[Why? Remember Lemma \ref{3.4x}(1) and the clause $(\gamma)$ of the
definition of $P^*_{\beta+1}$.] By Lemma \ref{3.4x}, we find a $(\pi^*(
\beta)+1)$--component $q'=\langle\zeta',\tau',n',u',\langle q_\vare', 
v_\vare':\vare<\theta\rangle\rangle$ of $p$ such that $\pi(i^*)\in u^{q'}$
and  
\[(\forall j\in u^{q'})(\forall\gamma\in (\pi^*(\beta),\rht(p)))(h^p(
\pi(i^*),\gamma)\geq\theta\ \Rightarrow\ h^p(j,\gamma)\geq\theta).\]
We claim that then 
\begin{enumerate}
\item[$(\boxtimes)$] \quad $(\forall j\in U[p,Z_{\rho_0}]\cap u^q)(\pi(j)\in  
u^{q'}\cap U[p,Z_{\rho_1}])$. 
\end{enumerate}
Why? Fix $j\in U[p,Z_{\rho_0}]\cap u^q$. Let $r,r'$ be components of $p$
such that $r\leq_{\rm pr}p$, $r'\leq_{\rm pr}p$, $\rht(r)=\beta+1$,
$\rht(r')=\pi^*(\beta)+1$ (so $r$ and $q$, and $r',q'$, are isomorphic). The
sets $Z_{\rho_0}\cap (\beta+1)$ and $Z_{\rho_1}\cap (\pi^*(\beta)+1)$ are
$p$--closed, and they have the same values of $\Upsilon$, and therefore
$U[p,Z_{\rho_0}\cap(\beta+1)]$ and $U[p,Z_{\rho_1}\cap(\pi^*(\beta)+1)]$ are
(order) isomorphic. Also, these two sets are included in $u^r$ and $u^{r'}$,
respectively. So looking back at our $j$, we may successively choose $j_0\in 
u^r\cap U[p,Z_{\rho_0}\cap (\beta+1)]$, $j_1\in u^{r'}\cap U[p,Z_{\rho_1}
\cap (\pi^*(\beta)+1)]$, and $j^*\in u^q$ such that
\begin{itemize}
\item $(\forall\gamma\leq\beta)(h^q(j,\gamma)=h^r(j_0,\gamma))$,
\item $(\forall\ell'\leq\ell)(h^r(j_0,\alpha_{\ell'})=h^{r'}(j_1,\pi^*(
\alpha_{\ell'})))$, and 
\item $(\forall\gamma\leq\pi^*(\beta))(h^{r'}(j,\gamma)=h^{q'}(j^*,
\gamma))$.
\end{itemize}
Then we have
\[(\forall\ell'\leq\ell)(h^q(j,\alpha_{\ell'})=h^{q'}(j^*,\pi^*(
\alpha_{\ell'}))\quad\mbox{and}\quad(\forall\gamma\in\pi^*(\beta)\setminus
Z_{\rho_1})(h^{q'}(j^*,\gamma)\geq\theta).\]
To conclude $(\boxtimes)$ it is enough to show that $\pi(j)=j^*$. If this
equality fails, then there is $\gamma<\rht(p)$ such that $\theta\neq h^p(
\pi(j),\gamma)\neq h^p(j^*,\gamma)\neq\theta$. If $\gamma\leq\pi^*(\beta)$,
then necessarily $\gamma\in Z_{\rho_1}$, and this is impossible (remember
$h^p(j,\alpha_{\ell'})=h^p(\pi(j),\pi^*(\alpha_{\ell'}))$ for $\ell'\leq
\ell$). So $\gamma>\pi^*(\beta)$. If $h^p(\pi(j),\gamma)=\theta+1$, then
$h^p(j^*,\gamma)<\theta$ and (by the choice of $q'$) $h^p(\pi(i^*),\gamma)<
\theta$. Then $\gamma\in Z_{\rho_1}$ and $h^p(i^*,(\pi^*)^{-1}(\gamma))<
\theta$, and also $h^p(i^*,(\pi^*)^{-1}(\gamma))=h^p(j,(\pi^*)^{-1}(\gamma))
=\theta+1$ (by the choice of $i^*$), a contradiction. Thus necessarily
$h^p(\pi(j),\gamma)<\theta$ (so $\gamma\in Z_{\rho_1}$) and therefore 
\[\theta>h^p(j,(\pi^*)^{-1}(\gamma))=h^p(i^*,(\pi^*)^{-1}(\gamma))=h^p(
\pi(i^*),\gamma)= h^p(j^*,\gamma)\]
(as the last is not $\theta$), again a contradiction. Thus the statement in
$(\boxtimes)$ is proven.    

Now we may finish considering the current case. By the definition of the
function $\Upsilon$ (and by the choice of $\rho_0,\rho_1$) we have 
\[\zeta'=\zeta_\ell,\quad\tau'=\tau_\ell,\quad n'=n_\ell,\quad\mbox{and }\
\pi[v_\xi]=v'_\xi\ \mbox{ for }\xi<\theta\]
(and $\pi\rest v_\xi$ is order--preserving). Therefore 
\[G^{\rho_0}_{\rho_1}(f)(\tau^*(x_i:i\in v_\xi))=f(\tau'(x_i:i\in v'_\xi))
\qquad\mbox{ (for every $\xi<\theta$).}\]
By the inductive hypothesis, $G^{\rho_0}_{\rho_1}(f)\rest u^{q_\xi} \in
F^{q_\xi}$ (for $\xi<\theta$), so as $f\in F^p$ (and hence $f\rest u^{q'}\in
F^{q'}$) we may conclude now that $G^{\rho_0}_{\rho_1}(f)\rest u^q\in F^q$.
\medskip

\noindent{\em Case 3:}\quad $\beta\in Z_{\rho_1}\setminus Z_{\rho_0}$\\
Similar.
\medskip

\noindent{\em Case 3:}\quad $\beta\in Z_{\rho_0}\cap Z_{\rho_1}$\\
If $U[p,Z_{\rho_0}]\cap u^q=\emptyset=U[p,Z_{\rho_1}]\cap u^q$, then
$G^{\rho_0}_{\rho_1}(f)\rest u^q\equiv 0$ and we are easily done. If one of
the intersections is non-empty, then we may follow exactly as in the
respective case (2 or 3).
\end{proof}

Now we may conclude the proof of the theorem. Since 
\[\bB^p\models\tau(x_i:i\in w_{\rho_0})<\tau(x_i:i\in w_{\rho_1}),\] 
we find $f\in F^p$ such that $f(\tau(x_i:i\in w_{\rho_0}))=0$ and
$f(\tau(x_i:i\in w_{\rho_1}))=1$. It should be clear from the definition of
the function $G^{\rho_0}_{\rho_1}(f)$ (and the choice of $\rho_0,\rho_1$)
that  
\[G^{\rho_0}_{\rho_1}(f)(\tau(x_i:i\in w_{\rho_0}))=1\quad\mbox{and}\quad
G^{\rho_0}_{\rho_1}(f)(\tau(x_i:i\in w_{\rho_1}))=0.\]
But it follows from Claim \ref{3.5.2} that $G^{\rho_0}_{\rho_1}(f)\in F^p$,
a contradiction.
\end{proof}

\begin{conclusion}
It is consistent that for some uncountable cardinal $\theta$ there is a
Boolean algebra $\bB$ of size $(2^\theta)^+$ such that 
\[\Dep(\bB)=\theta\quad\mbox{ but }\quad (\omega,(2^\theta)^+)\notin
\Dep_{\rm Sr}(\bB).\] 
\end{conclusion}

\begin{problem}
Assume $\theta<\lambda=\lambda^{<\lambda}$ are regular cardinals. Does there
exist a Boolean algebra $\bB$ such that $|\bB|=\lambda^+$ and for every
subalgebra $\bB'\subseteq\bB$ of size $\lambda^+$ we have $\Dep(\bB')=\theta$?
\end{problem}

\shlhetal
\end{document}